\documentclass[12pt]{article}
\usepackage[a4paper,top=2cm,bottom=2cm,left=3cm,right=3cm,marginparwidth=1.75cm]{geometry}
\usepackage{amsmath, amssymb}
\usepackage{graphicx, graphics}
\usepackage{amsthm}
\usepackage{hyperref}

\newcommand{\inst}[1]{}
\newcommand{\orcidID}[1]{}
\newcommand{\titlerunning}[1]{}
\newcommand{\authorrunning}[1]{}
\newcommand{\institute}[1]{}
\newcommand{\keywords}[1]{\noindent{\bf Keywords.} #1}
\newtheorem{lemma}{Lemma}
\newtheorem{theorem}{Theorem}
\theoremstyle{definition}
\newtheorem*{definition}{Definition}

\makeatletter
\def\blfootnote{\gdef\@thefnmark{}\@footnotetext}
\makeatother
\def\toolittle{\vspace{6ex}}  
\def\toomessy{\thispagestyle{empty}
  \blfootnote{Presented at Applied Geometric Algebras in Computer Science and Engineering (AGACSE) Conference 2024.}}

\usepackage{mdwmath}
\usepackage{mdwtab}
\usepackage{array}
\usepackage{epsfig}
\usepackage{booktabs}
\usepackage[caption=false]{subfig}
\usepackage{times}
\usepackage{float}

\newcommand\Z{{\mathbb Z}} 
\newcommand\N{{\mathbb N}} 
\newcommand\R{{\mathbb R}}

\newcommand\HB{{\mathbb H}} 
\newcommand\OB{{\mathbb O}} 
\newcommand\SB{{\mathbb S}}
\newcommand{\PB}[1]{{\mathbb P}_{#1}}
\newcommand{\Pf}{{\rm pf}} 
\newcommand{\Cl}{{\rm Cl}} 
\newcommand{\GA}[1]{{\mathcal G}_{#1}}
\newcommand{\GT}{{\rm G2}}
\newcommand{\SL}{{\rm SL}}
\newcommand{\PSL}{{\rm PSL}}   
\newcommand{\Spin}[1]{{\rm Spin#1}} 
\newcommand{\Pin}[1]{{\rm Pin#1}} 
\newcommand{\Comb}{\mathcal{C}}

\newcommand{\wedgep}{\mathbin{\land}}
\newcommand{\crossp}{\mathbin{\times}}
\newcommand{\dotp}{\mathbin{\cdot}}

\newcommand{\degs}{^\circ}
\newcommand{\va}{{\bf a}} 
\newcommand{\vb}{{\bf b}} 
 
\newcommand{\vv}[1]{{\bf #1}} 
\newcommand{\e}[1]{{\rm e}_{#1}}
\newcommand{\p}[1]{{\rm p}_{#1}} 
 
\renewcommand{\u}[1]{{\rm u}_{#1}}

\newcommand{\eqb}{\vspace{0pt}\[}
\newcommand{\eqe}{\]}
\newcommand{\eqnb}{\vspace{0pt}\begin{equation}}   
\newcommand{\eqne}{\end{equation}}   
\newcommand{\arb}{\vspace{0pt}\[\begin{array}{ll}}
\newcommand{\are}{\end{array}\]}
\newcommand{\arnb}{\vspace{0pt}\begin{equation}\begin{aligned}}
\newcommand{\arne}{\end{aligned}\end{equation}\vspace{0pt}}
\newcommand{\tabh}{\begin{table}[ht]}
\newcommand{\tabb}{\begin{table}}
\newcommand{\tabe}{\end{table}}

\newcommand{\hod}{^*}
\newcommand{\ec}{\text{,}}
\newcommand{\es}{\text{.}}

\newcommand\quid{\hskip 0.6em}
\newcommand\sm{\footnotesize -}
\newcommand\sn{\footnotesize}
\renewcommand\sm{\tiny\;-}
\renewcommand\sn{\tiny\;}
\newcommand\HH[1]{H\textsubscript{#1}}
\newcommand\XX[1]{X\textsubscript{#1}}
\newcommand\YY[1]{Y\textsubscript{#1}}
\newcommand\CJ{C{+}J} \newcommand\CJB{(C{+}J)\;}
\newcommand\LE{L{\sm}E} \newcommand\LEB{(L{\sm}E)\;}
\newcommand\DK{D{+}K} \newcommand\DKB{(D{+}K)\;}
\newcommand\FM{F{+}M} \newcommand\FMB{(F{+}M)\;}
\newcommand\AH{A{+}H} \newcommand\AHB{(A{+}H)\;}
\newcommand\GN{G{+}N} \newcommand\GNB{(G{+}N)\;}
\newcommand\BI{B{+}I} \newcommand\BIB{(B{+}I)\;}

\begin{document}
\title{Construction of exceptional Lie algebra \texorpdfstring{$\GT$}{G2} and
       non-associative algebras using Clifford algebra}
\titlerunning{Construction of $\GT$ and non-associative algebras using Clifford algebra}
%
\author{G.\:P.\:Wilmot\inst{}\orcidID{0009-0003-7451-3933}}
\authorrunning{G.\:P.\:Wilmot}
\institute{University of Adelaide, Adelaide, South Australia, 5005, Australia
\email{greg.wilmot@adelaide.edu.au}\\
\url{https://www.adelaide.edu.au}}
\maketitle              
\toomessy
\begin{abstract}
This article uses Clifford algebra of definite signature to derive octonions and the Lie exceptional algebra \texorpdfstring{$\GT$}{G2} from calibrations using \texorpdfstring{$\Pin{(7)}$}{Pin(7)}. This is simpler than the usual exterior algebra derivation and uncovers a subalgebra of \texorpdfstring{$\Spin(7)$}{Spin(7)} that enables \texorpdfstring{$\GT$}{G2} and an invertible element used to classify six other algebras which are found to be related to the symmetries of \texorpdfstring{$\GT$}{G2} in a way that breaks the symmetry of octonions. The 4-form calibration terms of \texorpdfstring{$\Spin{(7)}$}{Spin(7)} are related to an ideal with three idempotents and provides a direct construction of \texorpdfstring{$\GT$}{G2} for each of the {\rm 480} representations of the octonions. Clifford algebra thus provides a new construction of \texorpdfstring{$\GT$}{G2} without using the Lie bracket. This result is extended to 15 dimensions generating another 100 algebras as well as the sedenions.

\end{abstract}

\section{Introduction}
Clifford algebra, $\Cl(p,q)$, is related to geometry and extends the field of differential geometry when limited to definite signature, either $\Cl(n,0)$ or $\Cl(0,n)$. This article concentrates on the former positive signature, denoted $\GA{n}$, which expands onto the exterior algebra using the Pfaffian. Clifford algebra is defined as a quadratic form over the exterior space defined by Grassmann and a null quadratic form recovers the exterior product consisting of $n$-forms. The $n$-forms of differential geometry, $dx$ and $dy$, are independent by definition only, which makes them exterior. In this respect differential geometry is a subset of $\GA{n}$ and this allows the full power of Clifford algebra to be applied to problems in differential geometry.

Another reason to specify geometric algebra as the Clifford algebra of definite signature is that this provides a one-to-one relationship with simplices. An $n$-simplex is defined by Spivak as \cite{Spivak} 
\eqb \Delta_n = \left\{\vv{x}\in\R^{n+1}: 0\le x_i\le1\text{ and }\sum_{i=1}^{n+1}x_i=1\right\}\es\eqe
The edges of the $n$-simplex are the simplest way of connecting $n+1$ points in an $n$ dimensional space so the corresponding geometric algebra is the algebra of the simplest geometry. This is formalised using the Pfaffian, which is a construction of the connections between vertices of the simplex and provides the rigour to realise the Spin and Pin groups. Marcel Riesz, \cite{Riesz}, exposed the relationship between the geometric and exterior algebras by deriving $\GA{n}$ from the external algebra and vise versa, as constructed by P.\;Lounesto in \cite{Lounesto}. The Pfaffian expansion of $\GA{}$ formalises the mapping into the Grassmann algebras as the fundamental structure equation \cite{Cannello,Wilmot1,Wilmot2},
\eqnb \va_1 \va_2 \va_3 \dots \va_k = \sum^{[\frac{k}2]}_{i=0} \sum_{\mu\in\Comb} (-1)^\sigma \,\Pf(\va_{\mu_1}\dotp
\va_{\mu_2},\dots,\va_{\mu_{2i-1}}\dotp\va_{\mu_{2i}}) \va_{\mu_{2i+1}}\wedgep\dots\wedgep\va_{\mu_k}\ec\label{eqn:fun}\eqne
where $\Pf(A)$ is the Pfaffian of $A$ and $\Comb = \binom{k}{2i}$  provides combinations, $\mu$, of $k$ indices divided into $2i$ and $k-2i$ parts and $\sigma$ is the parity of the combination. Vectors, $\va_{k}$, in $n$ dimensions are linear combinations of basis elements, $\e{i}$, $i\in\N_1^n = \{1,2,3,\dots,n\}$, using the set notation of Porteous\cite{Porteous}.

The wedge product of Grassmann algebra, $\wedgep$, can only be represented in geometric algebra as a selection of the exterior part of the geometric product, which, in general, is a semi-graded product as shown in \cite{Wilmot1}. Equation (\ref{eqn:fun}) is derived from the basis product, $\e{i} \e{j} = \delta_{ij} + \e{i}\wedgep \e{j}$, $i,j\in\N_1^n$.  The number of independent vectors, $k$, in (\ref{eqn:fun}) is limited by uniqueness as described after introducing the Spin group.
\begin{definition} The Spin and Pin groups are even and odd products of unit length vectors in (\ref{eqn:fun}), respectively \cite{Harvey}. The cofactor expansion of the Pfaffian and determinant guarantees associativity of (\ref{eqn:fun}) which is thus invertible for non-zero vectors. The number of independent pairs of vectors correlates to the edges of the labelled $(n-1)$-simplex, for an $n$ dimensional space.\end{definition} 

\begin{figure}[ht] \begin{center}
  \subfloat[\centering Geometric algebra $\GA{3}$]{{
  \includegraphics{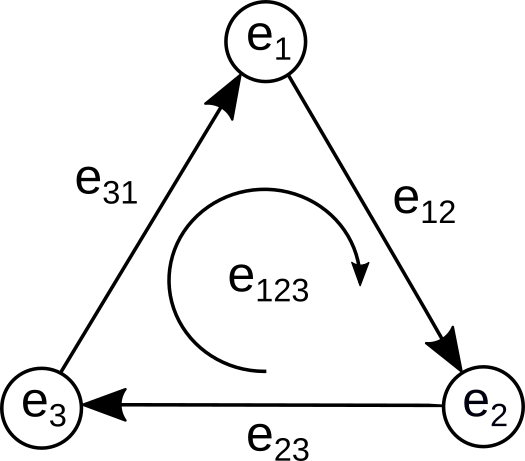} \label{fig:2sima}}}
  \qquad\qquad
  \subfloat[\centering Cayley-Dickson algebra $\HB$]{{
  \includegraphics{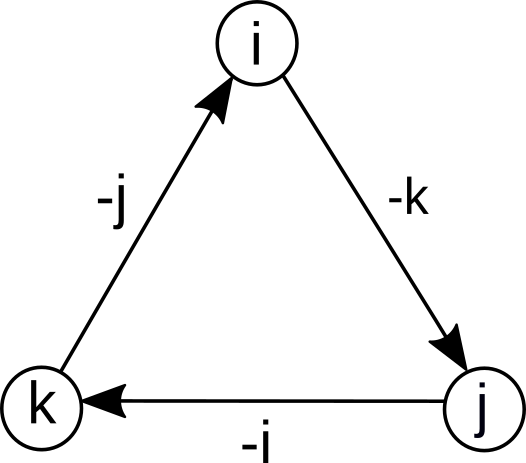} \label{fig:2simb}}}
  \caption{\bf The 2-simplex correspondences} \label{fig:2sim}
\end{center}
\end{figure}

Labelling simplices, using the notation derived from differential geometry, then translating the vertices into separate dimensions, the edges to areas, faces to volumes, etc provides a relationship to $\GA{}$, as shown for $\GA{3}$ in Figure \ref{fig:2sima}. The Pfaffian provides the mechanics connecting the two and Pascal's triangle, known to relate to simplices, provides the cardinality of the Grassmannians and for rotations this is $\binom{n}2$, which gives the number of unique pairs of vectors in (\ref{eqn:fun}). For $\GA{3}$ this is $3$ and Figure \ref{fig:2sima} shows these as the three basis edges of the $2$-simplex which forms a $3$-cycle, $\e{12}\e{23}\e{31} = 1$, related to quaternions. The fourth row of Pascal's triangle is $\{1,3,3,1\}$ which gives one scalar, three vectors, the three 2-forms just provided and one pseudoscalar 3-form. The pseudoscalar is related to the $3$-cycle via the cross product, $\va\crossp\vb = \e{321}\va\wedgep\vb$, shown here for the right-hand screw rule. The pseudoscalar in $\GA{3}$ is imaginary, in that it commutes with $\GA{3}$ and has negative square. But denoting it as $i$ is inadequate in the cross product because the parity would be ambiguous. Also $\GA{}$ has imaginary pseudoscalars in many other dimensions, such as $\GA{7}$ and $\GA{15}$, both referenced later.

General associativity of (\ref{eqn:fun}), $\va_1 \va_2\dots\va_{n} = (\va_1\dots\va_{r})(\va_{r+1}\dots\va_{n})$, $r\in\N_1^n$ \cite{Wilmot1}, allows pairs of vectors to be isolated, allowing the Spin group to be redefined as 
\eqnb \begin{array}{c} R = \prod_{i,j \in\Comb} R_{ij} (\theta_{ij})\text{ where }
R_{ij}(\theta) = \cos(\theta) +\e{ij} \sin(\theta)\ec
\end{array}  \label{eqn:rot}\eqne

where the combinatorial expansion $\Comb = \binom{n}2$ provides all combinations of pairs, $i,j\in\N_1^n$, $i < j$ and $\theta_{ij}$ are generalised Euler angles governed by the number of edges of the $n-1$ simplex. Hamilton defined the versor as the rotation in three dimensions using quaternions. Hence (\ref{eqn:fun}) for $k$ even can be called a generalised versor. With an odd number of vectors, (\ref{eqn:fun}) is extended to the Pin group which provides reflections remembering that any pair that overlaps with any pair of existing vectors just changes the corresponding Euler angle by $180\degs$. Both the Spin and Pin group act under conjugations as double coverings of the orthogonal group. For $R$ a rotation and $P$ a $p$-form in any $\GA{}$, the conjugations for rotations and reflections are

\arnb A' &= R A R^{-1} \quad\text{ and,} \\
      A'' &= (-1)^{(p+1)} P A' P^{-1}\es
\label{eqn:reflect} \arne

\begin{definition} $R_{ij} = R_{ij}(\frac{\pi}4) = \frac1{\sqrt2}(1 +\e{ij})$ defines the important set of conjugation rotations of $90\degs$ for basis elements. This sends any $\e{j}$ component in $A$ to $\e{i}$ and $\e{i}$ to $-\e{j}$. Using this swapping trick rather than applying the rotation $R_{ij}(\pi/2)$ avoids rounding errors being accumulated. \end{definition} 
\begin{definition} $R_{ijkl}$ is shorthand for all three $90\degs$ versor pairs; $R_{ij}R_{kl}$, $R_{ik}R_{jl}$ and $R_{il}R_{jk}$, $i,j,k,l$ all distinct. For example, $R_{ij} R_{kl} = \frac12(1 + \e{ij} + \e{kl} +\e{ijkl})$, and all three rotations involve $\e{ijkl}$. Of the 24 permutations of the indices only these three rotations need be considered in the conjugation. In dimensions greater than $5$ the product of three versors may involve 6-forms but in the following only 4-forms will be required for the general proofs below. \end{definition} 

Certain simplices are also used to generate Cayley-Dickson algebras, such as quaternions, $\HB$, as shown in Figure \ref{fig:2simb} using the usual formalism for the three rules $ij=k$, $jk=i$ and $ki=j$ for the vertices. Using the notation from $\GA{3}$, the edges can also be labelled but including vectors and rotations from the same algebra generates an inconsistency with opposite parities for vectors and rotations so in the following the usual notation of just labelling the vertices will be adopted. Cayley-Dickson algebras have negative signature so one of the vertices must be negated in order to represent an identity $3$-cycle. This opposite parity case, is isomorphic to reversed arrows and vertices labelled as $-i$, $-j$, $-k$ and is anti-isomorphic to Figure \ref{fig:2simb}. All combinations of arrow directions are isomorphic to one of these two algebras so only the two 3-cycle directions need be considered. 

This is the basis of the connection between $\GA{}$ and the Cayley-Dickson algebras. Certain $3$-forms can generate 3-cycles that realise a multiplication table that covers all products of the algebra apart from the squares. Taking the squares of the simplex vertices to be $-1$ and using the 3-cycle rules given by the parity of each basis 3-form it is only necessary to find $3$-forms with independent $3$-cycles that cover all products. This is found for the 6-simplex and 14-simplex, for example, which have much greater enumeration of parity than the 2 parities of quaternions and leads to new algebras related to split octonions and sedenions. In $\GA7$ they are easily visualised as symmetry operations of the 6-simplex and as partial symmetries of the Cartan root diagram for $\GT$. This leads to a subalgebra of $\Spin(7)$ that enables certain 3-forms, called calibrations, to generate octonions in differential geometry and provide a complete derivation of the Lie algebra $\GT$ from $\Spin(7)$. The automorphisms that keep the calibrations invariant are easily visualised in $\GA7$ while the new algebras are divided into six classes that are invariant to some of the $\GT$ operators. Some of these concepts carry over to $\GA{15}$ and sedenions.

The next section uncovers the enabling subalgebra from $\Spin(7)$ and uses this to classify the seven algebras related to calibrations. The section following this on octonion representations provides mappings for all $480$ octonions which demonstrates the conciseness of the algebra and provides completeness of the theorems. The section on sedenions extends the results to $\Spin(15)$ and provides examples of zero divisors that are associated with power associative algebras. The final section specifies the construct for any of the $480$ representations of the exceptional Lie algebra $\GT$ in and compares the Bryan and Cartan representations.

\section{Calibrations in GA(7)}
In $\GA{7}$, the Pascal's triangle row for the 6-simplex is $\{1, 7, 21, 35, 35, 21, 7, 1\}$ which corresponds to 7 basis vectors, 21 edges or Euler angles and 35 3-forms. Since the number of edges and faces are both divisible by 3 and 7, then 7 independent $3$-cycles can be selected to cover all edges once and all dimensions are covered 3 times. This defines a cross product in 7 dimensions, $\va\crossp\vb = \Phi_1\;\va\wedgep\vb$. Differential geometry represents this using the Fano plane with a certain arrangement of arrows whereas the flexibility of $\GA{}$ allows all possible arrangements of arrows to be explored. The Spin group is invaluable for this analysis which not only simplifies the calculations but finds new algebras and geometric ways of appreciating existing results. Such a cross product in $\GA{7}$ displayed with lexicographical ordering of the indices and with positive terms, called a primary 3-form, is

\eqb \Phi_1 = \e{123} +\e{145} +\e{167} +\e{246} +\e{257} +\e{347} +\e{356}\es \eqe

There are 30 ways to select primary cross product 3-forms which can be found by enumerating all combinations of 7 of the 35 triples of $\N_1^7$, eliminating those with duplicate edges. Alternatively, applying $90\degs$ rotations and reflections to $\Phi_1$ can be used, as will be done later. Table \ref{tab:class} lists the 30 primary 3-forms, called $\Phi_i$, $i\in\N_1^{30}$. The cross product is not reversible as it is not unique but does provide a consistent orthogonal vector for any 2-form, $\va\wedgep\vb$. For the moment consider the action of a single rotation, $R_{ij}$, which takes $\Phi_1$ to another 3-form in the primary list, apart from signs. Thinking geometrically, since each vertex of the 6-simplex has 3 touching primary faces then swapping the vertices of one face's edge will swap the same vertices of the other pairs of touching faces at each end of the reversed edge. But the pairs of faces at one end can not touch the faces at the remote end due to independence, guaranteeing each 3-form of the primary still has unique edges.

\begin{definition} The dual to any element of $\GA{7}$, $\Phi$, denoted $\Phi\hod$, uses the pseudoscalar
\eqb \Phi\hod = -\e{1234567}\;\Phi\es\eqe
This form is related to the Hodge star operator acting on $\Phi$ using the inner product of the Grassmann algebra. In geometric algebra this operation is just multiplication by the pseudoscalar and $\Phi_i\hod$, $i\in\N_1^{30}$, is a 4-form. 
\end{definition}

\begin{lemma}
The terms of $\{1 + \Phi_i\hod\}$ and $\{1 + \Phi_i + \Phi_i\hod+ \e{1234567}\}$, $i\in\N_1^{30}$, form commuting subalgebras of $\Spin(7)$ and $\Pin(7)$, respectively.
\end{lemma}
\begin{proof} By the definition of $\Phi_i$ a selected term, $\e{jkl}$, has each index occurring uniquely in a pair of other terms of $\Phi_i$ and no other shared index with $\e{jkl}$. Multiplying the selected term by the pair that shares $j$, $\e{jop} +\e{jqr}$ ($j, l, o, p, q, r$ unique), then $j$ will anti-commute twice (with $o$ and $p$ or $q$ and $r$) and $k$ and $l$ will each anti-commute three times. Thus giving even parity so each term of $\Phi_i$ commutes with $\Phi_i$ itself. Since $\e{1234567}$ commutes with $\GA{7}$ then $\Phi_i\hod$ also commutes with $\Phi_i$. Also the term $\e{ijk}$ when multipling $\Phi_i$ will contract with itself and form 4-forms with each of the other terms since there are pairs that share $i, j, k$ and the remaining indices are unique edges. But the result still commutes with $\Phi_i$. Now, again by the definition, any of the 28 3-forms not included in $\Phi_i$ must share an edge with one of the terms of $\Phi_i$ and will anti-commute with this term and its dual will also anti-commute. Hence the two sets are disjoint and $\Phi_i^2$ contains only terms from $\Phi_i\hod$ and a scalar and similarly for ${\Phi_i\hod}^2$, which are the second and first sets of the lemma, respectively.
\end{proof}

The dual of $\Phi_1$ is the 4-form 
\eqnb\Phi_1\hod = \e{1247} +\e{1256} +\e{1346} +\e{1357} +\e{2345} +e_{2367} +\e{4567}\es \label{eqn:form7}\eqne
We can now show that the $30$ primaries provide complete coverage of permutations of the vertices of the 6-simplex, $S_7$ which has order $7! = 30 \times 168$. The rotations provided by each $\Phi\hod$ have $7\times4! = 168$ possible rotations and this corresponds to the projective group of the cube $\PSL(2,\Z_7)$, isomorphic to $\SL(3,\Z_2)$, which has order $168$. Notice that the terms of each $\Phi\hod$ as $180\degs$ rotations leave $\Phi$ invariant because they commute. The permutations of the indices change the signs and follow the symmetric group $S_4$ that has order $24$ and act as the three pairs of rotations $R_{ijkl}$ for each term $\e{ijkl}$ along with $8$ sign variations. But the $30$ primaries were found by taking all $21$ rotations and the signs can be ignored if we consider all possible changes of arrow directions in the 6-simplex. Actually, the 2-simplex showed that the $8$ variation of arrows in a $3$-cycles corresponds to just two unique algebras so only changes of arrows in the Fano plane need be considered. Put another way, the $7!$ permutations of the vertices of the 6-simplex have $6\times5$ permutations of the $7$ selected independent faces that have no shared edges and the remaining $4!$ permutations are just rotations of some of the edges within these faces. Thus there are $(7!/168) \times 2^7 = 3,840$ possible label permutations and arrow combinations for the Fano plane, which are now considered.

\begin{figure}[ht] \begin{center}
  \begin{tabular}{c}
  \includegraphics{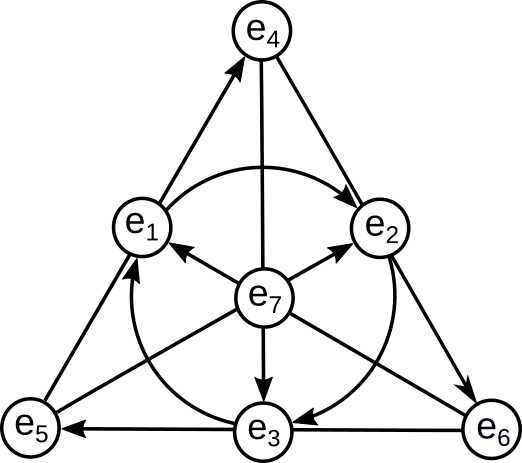}
  \end{tabular}
  \caption{\bf Primary cross product Fano plane} \label{fig:fano}
\end{center} \end{figure} 

The Fano plane in Figure \ref{fig:fano} shows a projection of the 6-simplex with arrows defined by $\Phi_1$. Like the circular arrow in $\GA{3}$, the arrows can define algebra elements and product rules but this diagram is not the usual representation that would define the product rules of octonions. Instead it denotes the natural cross product $\Phi_1$ and generates another non-associative algebra that will be exposed shortly. Changing the sign of $\e{246}$ in $\Phi_1$ means reversing the arrow on the right hand side of the diagram and this is now a representation of the octonion multiplication rule, denoted $\Phi_O$, called the associative calibration for $\OB$ and $\Phi\hod_O$ the coassociative calibration by Harvey and Lawson \cite{HarveyLawson}. There are 16 representations of octonions in each primary, which will be denoted $\Phi_{i,O}$, $i\in\N_1^{30}$. The well know exterior calibration product is \cite{HarveyLawson},
\eqnb \Phi_O\wedgep\Phi\hod_O = 7\e{1234567}\es\label{eqn:cal}\eqne 

In geometric algebra, the closure of the $\Phi_i\hod$ terms under multiplication means that $\Phi_i\Phi_i\hod$ can only consist of terms from $\Phi_i$ and the pseudoscalar. For calibrations the equivalent of (\ref{eqn:cal}) is,
\eqnb\Phi_O^2 + 7 = (-1)^{\sigma_j} 6\e{1234567} \Phi_O\ec \label{eqn:sqO}\eqne
where the $\sigma_j$ factor defines a sign parity for the octonions and classifies 8 $\OB^+$ and 8 $\OB^-$ algebras for $\sigma_j$ even or odd respectively. The inclusion of a metric in $\GA{}$ extends (\ref{eqn:cal}) to a form that is later shown to be related to a projection operator. This can be generalised further by considering not just the calibrations but all signed terms of the primary 3-forms. There are $2^7$ or 128 combinations of signs for the 7 terms of $\Phi_i$ and $\sigma_j$ is the number of minus signs in these combinations.
\begin{definition} The primary 3-forms of Table \ref{tab:class}, $\Phi_i$, for $i\in\N_1^{30}$ are extended to all combinations of signs of each term, $\Phi_{i,j}$, $j\in\N_1^{128}$. The sequence starts with all positive terms, $\Phi_i=\Phi_{i,1}$, followed by 7 single minus terms each in turn, then pairs, etc. Thus $\Phi_{1,5}$ denotes the $\OB^-$ algebra with $\e{246}$ negated in $\Phi_1$, as discussed above. It is only necessary to consider the first half terms because the second half is just the negation of the first half which swaps the $\OB^+$ and $\OB^-$ classes, $\Phi_{i,129-j} = -\Phi_{i,j}$ and in (\ref{eqn:sqO}), $\sigma_{129-j} = 7 - \sigma_j$ for $j\in\N_1^{64}$.
\end{definition}

This has lain the groundwork for classifying the seven non-associative algebras generated by the $30 \times 128$ representations provided by geometric algebra. A better way to represent $(\ref{eqn:sqO})$ is to define an invertible multivector, 
\eqb \rho_{i,j} = \frac14(3\e{1234567} -(-1)^\sigma\Phi_{i,j}),\eqe
where, if $\Phi_{i,j} = \Phi_O$ then $\rho_{i,j} = \rho_O$ and $\rho_O^2 = -1$, so $\rho_O^{-1} = -\rho_O$. With $\phi_{i,j}$, called the remainder, being one of the terms of $\Phi_i$ yet to be determined, then
\eqnb 2\e{1234567}(\rho_{i,j}^2 +1) = 
  \begin{cases} 0, \qquad\qquad\text{or} \\
       \pm\phi_{i,j} +(-1)^\sigma\Phi_{i,j}.
  \end{cases}
\label{eqn:rho}\eqne 
\begin{theorem}{\bf Classification Theorem}
For calibrations, $\rho_{i,j}^2 = -1$, otherwise $\pm\phi_{i,j}$ in $(\ref{eqn:rho})$ specifies via Table \ref{tab:class} that $\Phi_{i,j}$ generates a non-associative algebra with 4, 8, 10, 12, 14 or 16 non-associative products which shall be designated $\PB4, \PB8, \PB{10}, \PB{12}, \PB{14}$ or $\PB{16}$, respectively. Hence, the number of non-associative, unique triplets in the algebra $\PB{k}$ is $k$. Note that the octonion representations have 28 non-associative, unique triples and 7 associative triples. The sign parity, $\sigma$, is the number of minus signs in $\Phi_{i,j}$, compared to the primary $\Phi_i$. It is convenient to refer to both 0 and the $\phi_{i,j}$ as the remainder, as shown in Table \ref{tab:class} in the last column. Each of the signed 3-forms, $\Phi_{i,j}$, has a remainder that specifies a generated algebra $\OB$, if $0$, or $\PB{k}$ provided by the Classes column in Table \ref{tab:class} for each $\Phi_i$ term as remainder, in order. For example, the classes for $\Phi_1$ show that a remainder of $\e{123}$ means $\PB4$, the next two terms as remainders mean $\PB{12}$ and the remaining four mean $\PB{14}$. The third column in Table \ref{tab:class} shows the remainder for the primary so that $\phi_{1,1}=-\e{246}$ specifies $\PB{14}$. Note that $\Phi_{11}$ and $\Phi_{20}$ have no remainder so these primaries generate $\OB$.
\end{theorem}

\tabh\caption{Primary Table and Classification Map}
\label{tab:class} \small \centering
\begin{tabular}{|c|c|c|c|}  \hline   
${\bf i}$ &${\bf \Phi_i}$ &{\bf Classes} &{\bf $\phi_{i,1}$}\\ \hline
$1$&$\e{123}+\e{145}+\e{167}+\e{246}+\e{257}+\e{347}+\e{356}$&$(\PB4,2\PB{12},4\PB{14})$&$-\e{246}$\\\hline
$2$&$\e{123}+\e{145}+\e{167}+\e{247}+\e{256}+\e{346}+\e{357}$&$(\PB4,2\PB{12},4\PB{14})$&$-\e{357}$\\\hline
$3$&$\e{123}+\e{146}+\e{157}+\e{245}+\e{267}+\e{347}+\e{356}$&$(\PB4,2\PB{14},2\PB{12},2\PB{14})$&$-\e{157}$\\\hline
$4$&$\e{123}+\e{146}+\e{157}+\e{247}+\e{256}+\e{345}+\e{367}$&$(\PB4,4\PB{14},2\PB{12})$&$-\e{146}$\\\hline
$5$&$\e{123}+\e{147}+\e{156}+\e{245}+\e{267}+\e{346}+\e{357}$&$(\PB4,2\PB{14},2\PB{12},2\PB{14})$&$-\e{346}$\\\hline
$6$&$\e{123}+\e{147}+\e{156}+\e{246}+\e{257}+\e{345}+\e{367}$&$(\PB4,4\PB{14},2\PB{12})$&$-\e{257}$\\\hline
$7$&$\e{124}+\e{135}+\e{167}+\e{236}+\e{257}+\e{347}+\e{456}$&$(\PB8,2\PB{12},3\PB{14},\PB{10})$&$-\e{135}$\\\hline
$8$&$\e{124}+\e{135}+\e{167}+\e{237}+\e{256}+\e{346}+\e{457}$&$(\PB8,2\PB{12},3\PB{14},\PB{10})$&$-\e{124}$\\\hline
$9$&$\e{124}+\e{136}+\e{157}+\e{235}+\e{267}+\e{347}+\e{456}$&$(\PB8,2\PB{14},2\PB{12},\PB{14},\PB{10})$&$-\e{267}$\\\hline
$10$&$\e{124}+\e{136}+\e{157}+\e{237}+\e{256}+\e{345}+\e{467}$&$(\PB8,4\PB{14},\PB{12},\PB8)$&$-\e{237}$\\\hline
$11$&$\e{124}+\e{137}+\e{156}+\e{235}+\e{267}+\e{346}+\e{457}$&$(\PB8,2\PB{14},2\PB{12},\PB{14},\PB{10})$&$0$\\\hline
$12$&$\e{124}+\e{137}+\e{156}+\e{236}+\e{257}+\e{345}+\e{467}$&$(\PB8,4\PB{14},\PB{12},\PB8)$&$-\e{156}$\\\hline
$13$&$\e{125}+\e{134}+\e{167}+\e{236}+\e{247}+\e{357}+\e{456}$&$(2\PB{10},\PB{12},\PB{14},\PB{16},\PB{12},\PB{10})$&$-\e{247}$\\\hline
$14$&$\e{125}+\e{134}+\e{167}+\e{237}+\e{246}+\e{356}+\e{457}$&$(2\PB{10},\PB{12},\PB{14},\PB{16},\PB{12},\PB{10})$&$-\e{356}$\\\hline
$15$&$\e{125}+\e{136}+\e{147}+\e{234}+\e{267}+\e{357}+\e{456}$&$(\PB{10},\PB{14},\PB{16},\PB{10},2\PB{12},\PB{10})$&$-\e{136}$\\\hline
$16$&$\e{125}+\e{136}+\e{147}+\e{237}+\e{246}+\e{345}+\e{567}$&$(\PB{12},4\PB{14},\PB{12},\PB4)$&$-\e{136}$\\\hline
$17$&$\e{125}+\e{137}+\e{146}+\e{234}+\e{267}+\e{356}+\e{457}$&$(\PB{10},\PB{14},\PB{16},\PB{10},2\PB{12},\PB{10})$&$-\e{457}$\\\hline
$18$&$\e{125}+\e{137}+\e{146}+\e{236}+\e{247}+\e{345}+\e{567}$&$(\PB{12},4\PB{14},\PB{12},\PB4)$&$-\e{247}$\\\hline
$19$&$\e{126}+\e{134}+\e{157}+\e{235}+\e{247}+\e{367}+\e{456}$&$(\PB{12},\PB{10},\PB{14},\PB{12},\PB{16},2\PB{10})$&$-\e{134}$\\\hline
$20$&$\e{126}+\e{134}+\e{157}+\e{237}+\e{245}+\e{356}+\e{467}$&$(\PB{12},\PB{10},3\PB{14},\PB{12},\PB8)$&$0$\\\hline
$21$&$\e{126}+\e{135}+\e{147}+\e{234}+\e{257}+\e{367}+\e{456}$&$(2\PB{12},\PB{16},\PB{10},\PB{14},2\PB{10})$&$-\e{257}$\\\hline
$22$&$\e{126}+\e{135}+\e{147}+\e{237}+\e{245}+\e{346}+\e{567}$&$(\PB{12},4\PB{14},\PB{12},\PB4)$&$-\e{245}$\\\hline
$23$&$\e{126}+\e{137}+\e{145}+\e{234}+\e{257}+\e{356}+\e{467}$&$(\PB{12},2\PB{14},\PB{10},\PB{14},\PB{12},\PB8)$&$-\e{126}$\\\hline
$24$&$\e{126}+\e{137}+\e{145}+\e{235}+\e{247}+\e{346}+\e{567}$&$(\PB{12},4\PB{14},\PB{12},\PB4)$&$-\e{137}$\\\hline
$25$&$\e{127}+\e{134}+\e{156}+\e{235}+\e{246}+\e{367}+\e{457}$&$(\PB{12},\PB{10},\PB{14},\PB{12},\PB{16},2\PB{10})$&$-\e{235}$\\\hline
$26$&$\e{127}+\e{134}+\e{156}+\e{236}+\e{245}+\e{357}+\e{467}$&$(\PB{12},\PB{10},3\PB{14},\PB{12},\PB8)$&$-\e{467}$\\\hline
$27$&$\e{127}+\e{135}+\e{146}+\e{234}+\e{256}+\e{367}+\e{457}$&$(2\PB{12},\PB{16},\PB{10},\PB{14},2\PB{10})$&$-\e{146}$\\\hline
$28$&$\e{127}+\e{135}+\e{146}+\e{236}+\e{245}+\e{347}+\e{567}$&$(\PB{12},4\PB{14},\PB{12},\PB4)$&$-\e{135}$\\\hline
$29$&$\e{127}+\e{136}+\e{145}+\e{234}+\e{256}+\e{357}+\e{467}$&$(\PB{12},2\PB{14},\PB{10},\PB{14},\PB{12},\PB8)$&$-\e{357}$\\\hline
$30$&$\e{127}+\e{136}+\e{145}+\e{235}+\e{246}+\e{347}+\e{567}$&$(\PB{12},4\PB{14},\PB{12},\PB4)$&$-\e{246}$\\\hline
\end{tabular}\tabe

\begin{proof} The proof is split into two halves, the first looks at the eight signed combinations of all positive terms or a single negative term ($\Phi_{i,j}, j\in\N_1^8$) and it is seen that these cover all remainders. The second half extends these results to all 128 signs using reflections including negation. The proof starts with the observation that the remainder for each primary, $\Phi_{i,1}$, is the term that is negated to provide the first octonion for each primary, apart from $\Phi_{11}$ and $\Phi_{20}$. Considering $\Phi_1$, the permutations of the basis indices needed to proceed to the negation of the first term, $\e{123}$, can be found by forming the multiplication table for both $\Phi_{1,1}$ and $\Phi_{1,2}$, and searching all permutations of the first until an isomorphism to the second is found. Table \ref{tab:construct} shows the permutations from $\Phi_{1,1}$ to $\Phi_{1,2}$, etc up to $\Phi_{1,8}$, which moves the minus sign to each term in succession with the final permutation wrapping back to $\Phi_{1,1}$. The first to the second shows a single signed permutation because $\Phi$ is changing its signed parity, $\sigma$. Hence the first four rows have isomorphic permutations showing that they are the same algebra. The remaining rows have anti-isomorphic relationships that need further analysis. Note that the transformations also permute the remainders correctly, as defined by (\ref{eqn:rho}), and that row 5 is missing because $\Phi_{1,5}$ is a calibration which has no isomorphic or anti-isomorphic relationship to a $\PB{k}$ algebra so the construction skips such cases because no permutation can change the remainder to $0$. 

\tabb\caption{$\Phi_1$ Simple Sign Combinations and Transformation Constructions} 
\label{tab:construct} \small   \centering
\begin{tabular}{|c|c|c|c|c|c|} \hline
${\bf \;j\;}$
&${\bf \Phi_{1,j}}$
&{\bf Remainder}
&{\bf Permutation}
&{\bf Rotation}
&{\bf Reflection}\\ \hline
$1$ &$\e{123}+\e{145}+\e{167}+\e{246}$ &$-\e{246}$ &$(-1)(4 5)(6 7)$ &$R_{45}R_{67}$ &$\e{157}$\\
 &$+\e{257}+\e{347}+\e{356}$  &&&&\\ \hline
$2$ &$-\e{123}+\e{145}+\e{167}+\e{246}$  &$\e{257}$ &$(23)(45)$ &$R_{23}R_{45}$ &$\e{12467}$\\
 &$+\e{257}+\e{347}+\e{356}$  &&&&\\ \hline
$3$ &$\e{123}-\e{145}+\e{167}+\e{246}$  &$\e{347}$ &$(45)(67)$ &$R_{45}R_{67}$ &$\e{12346}$\\
 &$+\e{257}+\e{347}+\e{356}$  &&&&\\ \hline
$4$ &$\e{123}+\e{145}-\e{167}+\e{246}$  &$\e{356}$ &$(1463)(25)$ &$R_{14}R_{16}$ &$\e{127}$\\
 &$+\e{257}+\e{347}+\e{356}$  &&&$R_{13}R_{25}$&\\ \hline
$6$ &$\e{123}+\e{145}+\e{167}+\e{246}$  &$\e{123}$ &$(2\;{-}4)(3\;{-}5)$ &$R_{24}R_{35}$ &$1$\\
 &$-\e{257}+\e{347}+\e{356}$  &&&&\\ \hline
$7$ &$\e{123}+\e{145}+\e{167}+\e{246}$  &$\e{145}$ &$(-2)(-3)$ &$R_{46}R_{57}$ &$\e{145}$\\
 &$+\e{257}-\e{347}+\e{356}$  &&$(46)(57)$&&\\ \hline
$8$ &$\e{123}+\e{145}+\e{167}+\e{246}$ &$\e{167}$ &$(-12)(4576)$ &$R_{12}R_{45}$ &$\e{34}$\\
 &$+\e{257}+\e{347}-\e{356}$  &&&$R_{47}R_{46}$&\\ \hline
\end{tabular}\tabe

All permutations can be mapped to $90\degs$ rotations and reflections. Table \ref{tab:construct} shows permutations as cycles, with $(i, j, k)$ meaning $(i\rightarrow j, j\rightarrow k, k\rightarrow i)$ and rotations with reflections to change signs. Permutation indices with minus signs mean the term changes sign. Equivalently, 2-form rotations and reflections to change signs can be used to replace any permutation. For example, the permutation $(j,k) \rightarrow (k,j)$ is $R_{jk}$ and a reflection of basis k. Similarly, $(i,j,k,l)\rightarrow(j,k,l,i)$ can be expressed as rotations $R_{ij}R_{ik}R_{il}$ in left to right order followed by sign changes for j, k and l, which, from the second of (\ref{eqn:reflect}), is the reflection 4-form with j, k and l indices missing. Since these rotations involve the 4-forms from $\Phi_1\hod$, which commute with $\Phi_1$, then equation $(\ref{eqn:rho})$ is invariant apart from sign changes of $\Phi_{1,j}$ and $\phi_{i,k}$. Finally, reflections are used to correct the rotation signs as well as address permutation sign changes. Negative permutations change the term's sign if it contains the negated index.  For example, the permutation in the first row changes the sign of all terms in $\Phi_{1,1}$ that contain $e_1$ and swaps $e_4$ to $e_5$ and $e_6$ to $e_7$ plus in the opposite direction with a sign change. A no reflection case is shown for $\Phi_{1,6}$ which is 1 according to $(\ref{eqn:reflect})$.

Converting the transformations for each row into rotations and reflections for both $\Phi$ and the remainder shows that the form of (\ref{eqn:rho}) is invariant so the theorem still holds. This process has been verified for all primary simple sign combinations for the first 8 progressions for each $\Phi_i, i\in\N_1^{30}, i \ne 11, 20$, and also for $\Phi_{11,2}$ and $\Phi_{20,2}$. Note that this is not class equivalence because, for example, the first term of both $\Phi_{15}$ and $\Phi_{16}$ is $\e{125}$ but as remainders these correspond to different algebras, $\PB{10}$ and $\PB{12}$. The reflections can only change the sign of remainders in $(\ref{eqn:rho})$, so while the rotations may cycle through the $\PB{k}$ algebras, the reflections only move to equivalent signed $\PB{k}$ or $\OB$ representations. Reflections are now used to generate all 128 combinations for each $\Phi_i$. 

Table \ref{tab:basis} shows the reflections of the first 8 signed combinations of $\Phi_{1,j}, j\in\N_1^8$, with a single basis. Along with the $\Phi_{1,j}$ row, this gives 64 elements which are all unique as shall now be proved. By the definition of the primary 3-forms each $\Phi_{1,j}$ contains the reflection $e_i$ three times so each row reflection changes the sign of the remaining 4 terms. These terms are unique for each row and have different non-overlapping sign differences for each $\Phi_{1,j}$ so each term is unique. The first column has 4 negative terms for each row after $\Phi_{1,1}$ in the header. Other columns may have 3 or 5 negative terms or 1 in the header. Now consider the negated table starting with header $\Phi_{1,128}, \Phi_{1,127}, \dots, \Phi_{1,121}$. The first column has 7 negatives in the header followed by 3 negatives for each row. The remaining columns have 2 or 4 negative terms and 6 in the header. So only the cases of 3 and 4 negative terms need be considered since these may overlap with Table \ref{tab:basis}.

\tabb\caption{Single Basis Reflections}
\label{tab:basis}\small\centering
\begin{tabular}{|c|c|c|c|c|c|c|c|c|}  \hline   
{\bf Reflection} &${\bf\Phi_{1,1}}$ &${\bf\Phi_{1,2}}$ &${\bf\Phi_{1,3}}$ &${\bf\Phi_{1,4}}$ &${\bf\Phi_{1,5}}$ &${\bf\Phi_{1,6}}$ &${\bf\Phi_{1,7}}$ &${\bf\Phi_{1,8}}$\\ \hline
$\e1$ &$\Phi_{1,99}$ &$\Phi_{1,114}$ &$\Phi_{1,119}$ &$\Phi_{1,120}$ &$\Phi_{1,64}$ &$\Phi_{1,63}$ &$\Phi_{1,62}$ &$\Phi_{1,61}$\\ \hline
$\e2$ &$\Phi_{1,90}$ &$\Phi_{1,105}$ &$\Phi_{1,60}$ &$\Phi_{1,54}$ &$\Phi_{1,117}$ &$\Phi_{1,118}$ &$\Phi_{1,48}$ &$\Phi_{1,47}$\\ \hline
$\e3$ &$\Phi_{1,85}$ &$\Phi_{1,100}$ &$\Phi_{1,55}$ &$\Phi_{1,49}$ &$\Phi_{1,46}$ &$\Phi_{1,45}$ &$\Phi_{1,115}$ &$\Phi_{1,116}$\\ \hline
$\e4$ &$\Phi_{1,79}$ &$\Phi_{1,59}$ &$\Phi_{1,104}$ &$\Phi_{1,43}$ &$\Phi_{1,111}$ &$\Phi_{1,38}$ &$\Phi_{1,113}$ &$\Phi_{1,36}$\\ \hline
$\e5$ &$\Phi_{1,76}$ &$\Phi_{1,56}$ &$\Phi_{1,101}$ &$\Phi_{1,40}$ &$\Phi_{1,37}$ &$\Phi_{1,110}$ &$\Phi_{1,35}$ &$\Phi_{1,112}$\\ \hline
$\e6$ &$\Phi_{1,72}$ &$\Phi_{1,52}$ &$\Phi_{1,42}$ &$\Phi_{1,103}$ &$\Phi_{1,106}$ &$\Phi_{1,33}$ &$\Phi_{1,32}$ &$\Phi_{1,109}$\\ \hline
$\e7$ &$\Phi_{1,71}$ &$\Phi_{1,51}$ &$\Phi_{1,41}$ &$\Phi_{1,102}$ &$\Phi_{1,34}$ &$\Phi_{1,107}$ &$\Phi_{1,108}$ &$\Phi_{1,31}$\\ \hline
\end{tabular}\tabe

But these terms are mixed from the first column in one table and the other columns from the other table. The reflections in the first column start with all terms having the same sign so each reflection changes terms that do not contain the reflection. In order for the remaining columns in the alternate table to have the same number of terms then the single term from the header, $e_{ijk}$, that has a different sign must be negated. This means the reflection can not involve any of the basis indices, $i, j, k$, and the three terms that contain one of these basis indices and the reflection basis will not change sign. That leaves three remaining terms, each containing one of these indices that change sign and thus differ from $e_{ijk}$. This does not match the pattern from the first column where the three same sign terms share a basis index. Thus all terms of Table \ref{tab:basis} and its negation are unique and this covers all possible $2^7 = 128$ sign combinations of $\Phi_1$. This argument can be extended to be applicable to all primary 3-forms which interestingly all show the same structure. Finally, the count of associative or non-associative triples must remain intact to keep the same remainder in $(\ref{eqn:rho})$ due to $\Phi_{i,j}$ being squared. For $\Phi_1$ shown in Table \ref{tab:basis}, the first three columns are all $\PB{14}$, the next $\OB$, the sixth is $\PB{4}$ and the last two $\PB{12}$. This completes the proof of the Classification Theorem.
\end{proof}

The distribution from  column $\Phi_{1,5}$ of Table \ref{tab:basis} can be analysed to show five $\OB^+$ and three $\OB^-$ for indices $\N_1^{64}$. This is the same for all representations apart from $\Phi_{11}$ and $\Phi_{20}$ which have seven $\OB^-$ and one $\OB^+$. The negation of $\Phi_{1,j}$ for $j\in\N_1^{64}$ has $j\in\N_{65}^{128}$ and the change of parity swaps this distribution. This provides all 480 representations of the octonions. The representation of octonions selected by Baez, \cite{Baez}, corresponds to $\Phi_{11,1}$ which has the less common distribution but the advantage of having all positive terms.

Table \ref{tab:class} shows that the primary 3-forms can be divided up into groups of six having the same pseudo-octonion distribution. All $\Phi_i$ rows contain $\PB{12}$ and $\PB{14}$. For each group of six rows, the following extra classes can be found, $\{\PB4\}$,$\{\PB8, \PB{10}\}$, $\{\PB4, \PB{10}, \PB{16}\}$, $\{\PB4, \PB8, \PB{10}, \PB{16}\}$, $\{\PB4, \PB8, \PB{10}, \PB{16}\}$. Thus applying rotations from the terms of $\Phi_i\hod$ to $\Phi_{i,j}$, not including the octonions, will only generate classes within these groups. 

The algebra of the terms $\Phi_i\hod$ can be called the enabling algebra of the exceptional Lie algebra $\GT$. This is because they commute with $\Phi_O$ so are related to automorphisms of octonions. This will be explored in the section that constructs $\GT$ from $\GA{7}$. This section also uncovers partial $\GT$ symmetry for the $\PB{k}$ algebras, $k\in \{4,8,10,12,14,16\}$ because the remainder in (\ref{eqn:rho}) also commutes with $\Phi_{i,j}$. For this reason these algebras will be called pseudo-octonion algebras.

The reason why $\rho_{1,5}$ is invertible is because the dual of the 3-form can be expressed as $\Phi_{1,5}\hod = 1 -8\Sigma$ where $\Sigma = \frac18(1-\e{1247})(1+\e{1357})(1-\e{4567})$. These are 3 independent idempotents defining the projection operator $\Sigma^2 = \Sigma$ which means 
\arb(\rho_{1,5}\hod)^2 &= (\frac14(3 +\Phi_{1,5}\hod))^2\\
                  &= (1 -2\Sigma)^2\\
                  &= 1\es\are

Hence $\rho_O\hod=(1-2\Sigma)$ is a projection operator and this construction of $\Sigma$ works for 28 of the 35 combinations of 3 terms from $\Phi_{i,O}\hod$ for all 30 octonion representations. It does not work for the pseudo-octonion algebra representations nor for positions $(1,2,7)$, $(1,3,6)$, $(1,4,5)$, $(2,3,5)$, $(2,4,6)$, $(3,4,7)$ or $(5,6,7)$. For another example, the first three terms also define $\Sigma$, $\Phi_{1,5}\hod = 1 -(1-\e{1247})(1-\e{1256})(1-\e{1346})$. Hence the forms $\rho_O$, even though invertible, can not be expressed as a product of vectors and hence are not part of $\Pin(7)$.

\section{Octonion Representations}
This section validates the Classification Theorem by providing transformations from the first primary to the others and the reflections for all octonions since the proof concentrated on the first primary. This section also applies to the automorphisms of the octonions, the exceptional Lie algebra $\GT$. The same transformations for the octonions are applicable to $\GT$ to generate the 480 representations of $\GT$. 

\tabh\caption{Transformations from $\Phi_1$ for the Primary 3-forms}
\label{tab:transform} \small    \centering
\begin{tabular}{|c|c|c|c|c|c|} \hline   
{\bf Result} &{\bf Rotation} &{\bf Reflection} &{\bf Result} &{\bf Rotation} &{\bf Reflection}\\ \hline
$\Phi_1$ &$1$ &$1$ &$\Phi_{16}$ &$R_{35}R_{46}R_{12}R_{34}$ &$\e6$\\ \hline
$\Phi_2$ &$R_{23}R_{45}R_{67}$ &$\e{246}$ &$\Phi_{17}$ &$R_{35}R_{47}R_{24}$ &$\e2$\\
  &&&& $R_{23}R_{56}R_{57}$ &\\ \hline
$\Phi_3$ &$R_{12}R_{45}R_{67}$ &$\e{256}$ &$\Phi_{18}$ &$R_{35}R_{47}R_{46}$ &$\e{145}$\\ \hline
$\Phi_4$ &$R_{12}R_{57}$ &$\e{345}$ &$\Phi_{19}$ &$R_{12}R_{36}R_{34}$ &$\e{146}$\\ \hline
$\Phi_5$ &$R_{23}R_{57}$ &$\e{157}$ &$\Phi_{20,2}$ &$R_{36}R_{47}R_{13}R_{12} $&\\
  &&&& $R_{47}R_{46}$ &$\e{57}$\\ \hline
$\Phi_6$ &$R_{45}R_{46}R_{47}$ &$\e{123}$ &$\Phi_{21}$ &$R_{36}R_{57}R_{35}R_{47}$ &$\e3$\\ \hline
$\Phi_7$ &$R_{12}R_{34}R_{56}$ &$\e{146}$ &$\Phi_{22}$ &$R_{36}R_{34}R_{35}$ &$\e{256}$\\ \hline
$\Phi_8$ &$R_{34}R_{67}R_{13}$  & $1$ &$\Phi_{23}$ &$R_{36}R_{13}R_{12}$ &$\e7$\\
  &$R_{12}R_{45}R_{47}$ &  &&$R_{46}R_{47}$&\\ \hline
$\Phi_9$ &$R_{34}R_{56}R_{36}R_{57}$ &$\e7$ &$\Phi_{24}$ &$R_{36}R_{45}R_{12}R_{57}$ &$\e4$\\ \hline
$\Phi_{10}$ &$R_{34}R_{56}R_{57}$ &$\e{367}$ &$\Phi_{25}$ &$R_{37}R_{46}R_{34}R_{56}$ &$\e5$\\ \hline
$\Phi_{11,2}$ &$R_{34}R_{56}R_{67}$ &$\e{234567}$ &$\Phi_{26}$ &$R_{23}R_{25}R_{27}+ $ &$1$\\
  &$R_{12}R_{23}R_{47}R_{45}$ & & &$R_{34}R_{36}R_{37}R_{35}$ &\\ \hline
$\Phi_{12}$ &$R_{34}R_{57}R_{12}R_{35}$ &$\e3$ &$\Phi_{27}$ &$R_{12}R_{37}R_{35}$ &$\e{356}$\\ \hline
$\Phi_{13}$ &$R_{35}R_{67}$ &$\e{257}$ &$\Phi_{28}$ &$R_{37}R_{56}R_{12}R_{34}$ &$\e3$\\ \hline
$\Phi_{14}$ &$R_{23}R_{24}R_{25}$ &$\e{167}$ &$\Phi_{29}$ &$R_{37}R_{45}R_{23}R_{67}$ &$\e1$\\ \hline
$\Phi_{15}$ &$R_{35}R_{46}R_{47}$   &$\e7$ &$\Phi_{30}$ &$R_{37}$ &$\e{123}$\\
  &$R_{12}R_{34}R_{37}R_{36}$ &&&&\\ \hline
\end{tabular}\tabe
Firstly, consider transformations that convert $\Phi_1$ to $\Phi_{i,1}, i\in\N_1^{30}, i\ne11,20$, and $\Phi_1$ to $\Phi_{11,2}$ and $\Phi_{20,2}$. These are provided in Table \ref{tab:transform} which was generated by considering the permutations of the isomorphism between the multiplication tables generated by $\Phi_1$ and $\Phi_i$, etc. This process follows the same rotation and reflection scheme as for the signed variations of Table \ref{tab:construct} but here some extra rotations are needed because reflections alone do not necessarily return the generated $\Phi_i$ to a primary. These cases are shown in Table \ref{tab:transform} with either no reflection, a single basis or as a 6-form. 

\tabh\caption{Octonion Reflections within each Primary}
\label{tab:reflect}\small    \centering
\begin{tabular}{|c|c|} \hline   
 ${\bf Result}$ &{\bf Reflection for each Octonion} \\ \hline
$\Phi_{1,64}$ &$\e{14567}, \e{124567}, \e{123467}, \e{123456}, \e{23456}, \e{23467}, \e{24567}$ \\ \hline
$\Phi_{2,61}$ &$\e{14567}, \e{134567}, \e{123567}, \e{123457}, \e{23457}, \e{23567}, \e{34567}$ \\ \hline
$\Phi_{3,54}$ &$\e{24567}, \e{124567}, \e{123567}, \e{123457}, \e{23467}, \e{23456}, \e{14567}$ \\ \hline
$\Phi_{4,60}$ &$\e{24567}, \e{124567}, \e{123456}, \e{123467}, \e{23567}, \e{23457}, \e{14567}$ \\ \hline

$\Phi_{5,62}$ &$\e{14567}, \e{134567}, \e{123456}, \e{123467}, \e{23467}, \e{23456}, \e{34567}$ \\ \hline
$\Phi_{6,63}$ &$\e{14567}, \e{124567}, \e{123567}, \e{123457}, \e{23457}, \e{23567}, \e{24567}$ \\ \hline
$\Phi_{7,60}$ &$\e{23567}, \e{123567}, \e{123457}, \e{123456}, \e{24567}, \e{23467}, \e{13567}$ \\ \hline
$\Phi_{8,59}$ &$\e{23467}, \e{123467}, \e{123456}, \e{123457}, \e{34567}, \e{23567}, \e{13567}$ \\ \hline
$\Phi_{9,63}$ &$\e{13567}, \e{123567}, \e{124567}, \e{123467}, \e{23467}, \e{24567}, \e{23567}$ \\ \hline
$\Phi_{10,63}$ &$\e{13567}, \e{123567}, \e{123457}, \e{123467}, \e{23467}, \e{23457}, \e{23567}$ \\ \hline
$\Phi_{11,50}$ &$\e{23456}, \e{123456}, \e{123467}, \e{123457}, \e{3567}, \e{234567}, \e{134567}$ \\ \hline
$\Phi_{12,54}$ &$\e{23567}, \e{123567}, \e{124567}, \e{123456}, \e{23457}, \e{23467}, \e{13567}$ \\ \hline

$\Phi_{13,63}$ &$\e{13467}, \e{123467}, \e{124567}, \e{123457}, \e{23457}, \e{24567}, \e{23467}$ \\ \hline
$\Phi_{14,62}$ &$\e{13467}, \e{134567}, \e{123567}, \e{123456}, \e{23456}, \e{23567}, \e{34567}$ \\ \hline
$\Phi_{15,60}$ &$\e{23467}, \e{123467}, \e{123567}, \e{123456}, \e{24567}, \e{23457}, \e{13467}$ \\ \hline

$\Phi_{16,60}$ &$\e{23467}, \e{123467}, \e{123456}, \e{123567}, \e{24567}, \e{23457}, \e{13467}$ \\ \hline
$\Phi_{17,17}$ &$\e{13467}, \e{134567}, \e{124567}, \e{123457}, \e{23457}, \e{24567}, \e{34567}$ \\ \hline
$\Phi_{18,63}$ &$\e{13467}, \e{123467}, \e{124567}, \e{123457}, \e{23457}, \e{24567}, \e{23467}$ \\ \hline
$\Phi_{19,60}$ &$\e{23457}, \e{123457}, \e{123467}, \e{123456}, \e{24567}, \e{23567}, \e{13457}$ \\ \hline
$\Phi_{20,50}$ &$\e{23567}, \e{123567}, \e{123456}, \e{123467}, \e{3457}, \e{234567}, \e{134567}$ \\ \hline
$\Phi_{21,63}$ &$\e{13457}, \e{123457}, \e{124567}, \e{123567}, \e{23567}, \e{24567}, \e{23457}$ \\ \hline
$\Phi_{22,63}$ &$\e{13457}, \e{123457}, \e{124567}, \e{123456}, \e{23456}, \e{24567}, \e{23457}$ \\ \hline
$\Phi_{23,59}$ &$\e{23456}, \e{123456}, \e{123567}, \e{123467}, \e{34567}, \e{23457}, \e{13457}$ \\ \hline
$\Phi_{24,60}$ &$\e{23457}, \e{123457}, \e{123467}, \e{123567}, \e{24567}, \e{23456}, \e{13457}$ \\ \hline
$\Phi_{25,64}$ &$\e{13456}, \e{123456}, \e{123567}, \e{123457}, \e{23457}, \e{23567}, \e{23456}$ \\ \hline
$\Phi_{26,61}$ &$\e{13456}, \e{134567}, \e{124567}, \e{123467}, \e{23467}, \e{24567}, \e{34567}$ \\ \hline
$\Phi_{27,54}$ &$\e{23456}, \e{123456}, \e{124567}, \e{123467}, \e{23567}, \e{23457}, \e{13456}$ \\ \hline
$\Phi_{28,60}$ &$\e{23456}, \e{123456}, \e{123457}, \e{123567}, \e{24567}, \e{23467}, \e{13456}$ \\ \hline
$\Phi_{29,62}$ &$\e{13456}, \e{134567}, \e{123457}, \e{123567}, \e{23567}, \e{23457}, \e{34567}$ \\ \hline
$\Phi_{30,63}$ &$\e{13456}, \e{123456}, \e{124567}, \e{123467}, \e{23467}, \e{24567}, \e{23456}$ \\ \hline
\end{tabular}\tabe

Applying the Table \ref{tab:transform} transformation to $\Phi_{1,64}\in\Phi_{1,O}$ will produce octonion representations for each primary and it is a simple matter to use reflections to negate some terms to uncover the 16 octonion representations for each primary. Table \ref{tab:reflect} shows the result of applying each Table \ref{tab:transform} transformation to $\Phi_{1,64}$ then the reflections required to obtain the first half octonions in each primary. The second half is achieved by negation or reflecting with $\e{1234567}$. Note that the even reflections mark where parity changes occur. 

\toolittle
\section{Non-associative algebras and Zero Divisors}
Zero divisors are defined as two non-zero expressions that multiply to give zero. This is a property of power-associative algebras such as the sedenions, $\SB$. For algebras with products having unique triples, as given by the Fano plane, at a minimum the following must be considered as a definition of zero divisors, 
\eqb (a+b)(c+d) = 0\ec\text{ for }a,b,c,d\in\SB\ec\text{ all unique and pure or non-scalar.}\eqe
The octonions are non-associative and are normed with every expression invertible hence can not have any zero divisors. There are 28 non-associative triples of pure octonions leaving 7 associative products that are all quaternion-like algebras. These correspond to the lines (and circle) of the Fano plane. Octonions are generated from any three pure elements not all from the same quaternion-like algebras. Comparing this to the split octonions, ${\tilde\OB}$, which also have 28 non-associative pure triples but have four unitary elements, $\u{i}^2=1$, $i\in\N_1^4$, then this realises ideals of the algebra or projection operators $(\frac12(1\pm\u{i}))^2 = \frac12(1\pm\u{i})$ and the zero divisors $(1+\u{i})(1-\u{i}) = 0$.
Again they can be generated from three not completely associative pure basis elements with either 1, 2 or 3 unitary elements. All such choices are isomorphic \cite{Harvey} so four pure unitary elements can be selected to provide four ideals. Multiplying the zero divisors on each side or both sides by pure elements generates another 36 unique zero divisors. Of these, 12 do not contain $1$ or $-1$.

Hence the split octonions are power-associative even though they have the same number of non-associative triples as octonions. Octonions are alternate-associative so that swapping the terms in $(ab)c-a(bc)$ changes its sign. This is not true for power associative algebras. The extreme example is since algebras with negative signature basis have invertible elements, $(a+b)(a+b) \propto -1$, then multiplying this by $(c+d)$ just scales this term. But if this is a zero divisor then multiplying the right hand side first gives $(a+b)((a+b)(c+d)) = (a+b)0 = 0$. The pseudo-octonion algebras are power-associative with each having 12 zero divisors, for all $\PB{k}$, $k\in\{4,8,10,12, 14, 16\}$, even though they have different numbers of non-associative triples so they are analogous to split-octonions.

An example pseudo-octonion algebra is the representation $\Phi_{1,6}$ which has remainder $\e{123}$ and, from Table \ref{tab:class}, represents $\PB4$. The multiplication table for $\PB4$, as generated from the product rules dictated by $\Phi_{1,6}$, is shown in Table \ref{tab:multiple}. The diagonal terms are not specified by the rules derived from the Fano plane but to generate octonions and not split octonions, then negative squares are selected.

\tabb\caption{Multiplication Table for $\PB4$}
\label{tab:multiple}\small \centering
\begin{tabular}{|c|c|c|c|c|c|c|c|} \hline   
$\PB4^2$ &${\bf p_1}$ &${\bf p_2}$ &${\bf p_3}$ &${\bf p_4}$ &${\bf p_5}$ &${\bf p_6}$ &${\bf p_7}$\\ \hline
${\bf p_1}$ &$-1$ &$\p3$ &$-\p2$ &$\p5$ &$-\p4$ &$\p7$ &$-\p6$\\ \hline
${\bf p_2}$ &$-\p3$ &$-1$ &$\p1$ &$\p6$ &$-\p7$ &$-\p4$ &$\p5$\\ \hline
${\bf p_3}$ &$\p2$ &$-\p1$ &$-1$ &$\p7$ &$\p6$ &$-\p5$ &$-\p4$\\ \hline
${\bf p_4}$ &$-\p5$ &$-\p6$ &$-\p7$ &$-1$ &$\p1$ &$\p2$ &$\p3$\\ \hline
${\bf p_5}$ &$\p4$ &$\p7$ &$-\p6$ &$-\p1$ &$-1$ &$\p3$ &$-\p2$\\ \hline
${\bf p_6}$ &$-\p7$ &$\p4$ &$\p5$ &$-\p2$ &$-\p3$ &$-1$ &$\p1$\\ \hline
${\bf p_7}$ &$\p6$ &$-\p5$ &$\p4$ &$-\p3$ &$\p2$ &$-\p1$ &$-1$\\ \hline
\end{tabular}\tabe

The non-associative triples of this algebra refer to faces of the 6-simplex not included in the simplex faces
\eqb (\p4, \p5, \p6), (\p4, \p5, \p7), (\p4, \p6, \p7), (\p5, \p6, \p7)\es\eqe

Given the same representation these are a subset of the octonion nonassociate pure triples which is why the name pseudo-octonion is appropriate. The 12 zero divisors for Table \ref{tab:multiple} are

\eqb\begin{array}{c}
 (\p7 +\p1)(\p2 +\p4)\ec\ (\p6 +\p1)(\p2 +\p5)\ec\ (-\p5 +\p1)(\p2 +\p6)\ec\\
 (-\p4 +\p1)(\p2, \p7)\ec\ (-\p6 +\p1)(\p3 +\p4)\ec\ (\p7 +\p1)(\p3 +\p5)\ec\\
 (\p4 +\p1)(\p3 +\p6)\ec\ (-\p5 +\p1)(\p3 +\p7)\ec\ (\p5 +\p2)(\p3 +\p4)\ec\\
 (-\p4 +\p2)(\p3 +\p5)\ec\ (\p7 +\p2)(\p3 +\p6)\ec\ (-\p6 +\p2)(\p3 +\p7)\es
\end{array}\eqe

Each pseudo-octonion algebra has a different set of 12 zero-divisors but given these do not have unity elements then these power-associative algebras can be classed as forerunners of sedenions rather than being similar to split octonions. For each primary, Table 4 has the same structure for all six pseudo-octonion algebras as well as the octonions, it is only the signs that change. Hence, since all combinations of primary labels and signs for an anti-symmetric multiplication table have been covered, there are only 7 possible imaginary algebras with 3 generators.

Sedenions, $\SB$, are the next algebra after octonions in the Cayley-Dickson construction and are power-associative so are not normed. They have 252 non-associative triples and it is well known that they have 84 zero divisors \cite{Cawagas}. Geometric algebra also provides a construction for this algebra. Like $\GA{7}$, in 15 dimensions $\GA{15}$ also has a commuting imaginary pseudoscalar and a $\GT$-like enabling algebra. Using hexadecimal notation so that $\e{F}$ is the basis of the 15\textsuperscript{th} dimension, then the pseudoscalar is $\e{123456789ABCDEF}$, and a 3-form can be defined that covers all edges of the 14-simplex once, and all dimensions are listed 7 times, hence this defines a cross product in 15 dimensions with 35 terms, 
\eqb\begin{array}{cl}
\Phi =
 &\e{123} +\e{145} +\e{167} +\e{189} +\e{1AB} +\e{1CD} +\e{1EF} +\e{246} +\e{257} \\
 &+\e{28A} +\e{29B} +\e{2CE} +\e{2DF} +\e{347} +\e{356} +\e{38B} +\e{39A} +\e{3CF} \\
 &+\e{3DE} +\e{48C} +\e{49D} +\e{4AE} +\e{4BF} +\e{58D} +\e{59C} +\e{5AF} +\e{5BE} \\
 &+\e{68E} +\e{69F} +\e{6AC} +\e{6BD} +\e{78F} +\e{79E} +\e{7AD} +\e{7BC}\es
\end{array} \eqe

For each algebra generated from this using the triple product rules and each term being negated in turn then counting the number of non-associative products gives the following association
 \eqb  \begin{array}{c}
   210, 210, 210, 210, 210, 210, 210, 210, 224, 200, 224, 200, \\
   224, 200,  208, 208, 208, 208, 208, 208, 252, 228, 204, 180, \\
   236, 236, 188, 188, 220, 196, 220, 196, 204, 204, 204, 204\es
 \end{array} \eqe

This represents 12 unique, non-associative algebras. Notice that the 21st term, $\e{48C}$, when negated generates an algebra with 252 non-associative triples and this is isomorphic to sedenions, $\SB$. The term $-\e{48C}$ is analogous to $-\e{246}$ in $\GA{7}$. Applying all 105 $90\degs$ rotations repeatedly generates millions of primaries and the process above uncovers $100$ distinct pseudo-sedenion algebras with non-associative unique triples numbered in the hundreds.

Another result from geometric algebra is that for any calibration, $\Phi_O$, we can define $\Psi_O$ to be the 21 3-forms of $\GA{7}$ not included in the terms of $\Phi_O$, then with the correct signs we find
\eqb \Phi_O\hod\Psi_O = \Psi_O\es\eqe

For $\Phi_{1,5}$ we have
\eqb \begin{array}{cl} \Psi_O = &-\e{124} +\e{125} +\e{126} +\e{127} +\e{134} +\e{135} -\e{136} \\
                              &+\e{137} +\e{146} +\e{147} +\e{156} +\e{157} -\e{234} +\e{235} \\
                              &+\e{236} +\e{237} +\e{245} +\e{247} +\e{256} -\e{267} +\e{345} \\
                              &-\e{346} +\e{357} +\e{367} +\e{456} +\e{457} +\e{467} -\e{567}.
\end{array} \eqe
A similar result appears in $\GA{15}$.

\section{Construction of \texorpdfstring{$\GT$}{G2}}
This section covers the automorphisms of the octonions, known as the exceptional Lie algebra $\GT$. These are derived from the Lemma showing that the terms of $\Phi_i\hod$ commute with the primary and hence with the 16 representations of octonions, $\Phi_{i,O}$, for each $i\in\N_1^{30}$. The terms of each $\Phi_{i,O}\hod$, $\e{jklm}$ say, specify three rotations that define the enabling algebra for $\GT$ because they commute with $\Phi_{i,O}$ which corresponds to conjugate rotations of $90\degs$, twice, $R_{jk}(\frac\pi2) R_{lm}(\frac\pi2)=\e{jklm}$. Due to octonions having no remainder in the Classification theorem, any rotation results in an octonion representation, hence a pair of $90\degs$ rotations from the enabling algebra returns to the same primary, but not necessarily the starting point. Viewing the Fano plane in Figure \ref{fig:fano}, the arrows will have traded places, exposing the symmetry operations. For example $R_{12}$ swaps $\e1$ and $\e2$ and the edge $\e{246}$ becomes $\e{146}$ in the Fano plane but $R_{47}$ swaps $\e4$ and $\e7$ returning $\e{246}$ as a triple through the middle of the Fano plane. Remember that the arrows represent associative, quaternion-like products so swapping $\e1$ and $\e2$ turns $\e{246}$ into a non-associative product from a different primary representation. The other way of looking at this is that $R_{12}$ has the same action as $R_{47}$ and they both lead to the same representation. This leads to the following theorem.

\begin{theorem}{\bf Automorphism Theorem}
The enabling algebra of terms $\e{jklm}$, $j\ec k\ec l\ec m\in\N_1^7$ all unique, provides automorphisms of octonions within the equivalent primary form, as
\arb R_{jklm} \Phi_{i,O}  R_{mlkj} &= \Phi_{i,O}'\text{\rm\,, \quad where $\e{jklm}$ is a term of $\Phi_i\hod$, and}\\
   \hfill \Phi_{i,O}' &= \Phi_{i,O} \text{\rm, \quad if $\e{jklm}$ is a term of $\Phi_{i,O}\hod$}\es\are
\end{theorem}
\begin{proof}
The three rotations of $R_{jklm}$ are called normal, mixed and outer, for $R_{jk}R_{lm}$, $R_{jl}R_{km}$ and $R_{jm}R_{kl}$, respectively. The pairs of 2-forms for these three rotations can be called enabling rotations. The normal rotation expansion, this time including a parity, shows how important these are,
\arnb R_{jk} R_{lm} &= (\cos(\frac\pi2) +\e{jk} \sin(\frac\pi2)) 
                       (\cos(\frac\pi2) +(-1)^\sigma\e{lm} \sin(\frac\pi2)) \\
                   &= \frac12(1 +\e{jk} +(-1)^\sigma\e{lm} +(-1)^\sigma\e{jklm})\ec\label{eqn:auto}\arne
where $\sigma = \sigma(\Phi_{i,O}) + \sigma(\pm\e{jklm})+\sigma(R_{jklm})$ and $\sigma(x)$ is the parity of $x$. The last term is odd for mixed and even for normal and outer rotations and $\pm\e{jklm}$ is the term from $\Phi_{i,O}\hod$. 

This is divided into commuting terms $\alpha = \frac12(1\pm\e{jklm})$ and terms that rotate parts of $\Phi_{i,O}$, $\beta = \frac12(\e{jk}\pm\e{lm})$. These have the following rules which can easily be proved and have the same formulation for all enabling rotations
\eqb \alpha^2 = \alpha\text{,\qquad }\beta^2 = \alpha -1\text{\quad and \quad}\alpha\beta = \beta\alpha = 0\es\eqe
By the definition of $\Phi_i$, it consists of the dual of the rotation, $\e{jklm}\hod$, and the three enabling rotations with a shared index from $\Phi_i$, $\e{jklm}\hod$ called the cross terms. For example, two of the terms of $\Phi_i$ are $\pm\e{jkx}$ and $\pm\e{lmx}$, where $x$ is one of the dual's indices. This means the $90\degs$ rotations, $\e{jk}\Phi_i\e{kj}$ and $\e{lm}\Phi_i\e{ml}$, are the same because they act in the enabling pair in the same way, apart from signs. If this was not the case then the terms of $\beta$ would produce conflicting results not related to a primary. Hence we have two results from this, the first being the definition of $\Phi_{i,O}'$,
\eqb\Phi_{i,O}\beta = \beta\Phi_{i,O}'\text{\quad and \quad}\alpha\Phi_{i,O} = \alpha\Phi_{i,O}'\es\eqe
The second result follows because the cross terms will swap places and keep the same primary form. Because $\alpha = \frac12(1\pm\e{jklm})$, all cross terms cancel leaving $\frac12\e{1234567} \allowbreak(\pm\e{jklm} -1)$. This is the reason why $\rho_O^2 = -1$ since it was defined in such a way that each effective 4-form will add itself and $-1$ to the scaled result and the accumulated $\Phi_{i,O}\hod$ in (\ref{eqn:rho}) cancels by the definition of $\rho_i$. Hence $(\Phi_{i,O}')^2 = -1$.

So the rotations by $90\degs$ of $\Phi_{i,O}$ as a conjugation, with inverse $R_{ml} R_{kj}=\alpha-\beta$, are
\arb R_{jk} R_{lm} \Phi_{i,O} R_{ml} R_{kj} &= (\alpha+\beta)\Phi_{i,O}(\alpha-\beta) \\
                &= (\alpha^2\Phi_{i,O}+\beta\alpha\Phi_{i,O} -\Phi_{i,O}\alpha\beta -\beta\Phi_{i,O}\beta) \\
                &= \alpha\Phi_{i,O} -\beta^2\Phi_{i,O}' \\
                &= \alpha\Phi_{i,O} -(\alpha -1)\Phi_{i,O}' \\
                &= \Phi_{i,O}'\es\are
Now the parity for $\beta$ must be the same for $\Phi_{i,O}'$ and $\Phi_{i,O}$ and hence for the $\e{jklm}\hod$ term. Since this can be any term from $\Phi_{i,O}$ then the two are equal. The parity allows the cross terms to keep the same signs. This completes the proof and shows that the enabling 2-forms, $\frac12(\e{jk} \pm\e{lm})$, $\frac12(\e{jl} \mp\e{km})$ and $\frac12(\e{jm}\pm\e{kl})$, realise the action of the automorphisms $R_{jklm}$ and this is the basis of the definition of $\GT$ as the automorphisms of the octonions. 
\end{proof}
The 240 representations of the octonions generated from Tables \ref{tab:transform} and \ref{tab:reflect} can be applied to the $\GT$ algebra generated from $\Phi_{1,64}$ to give 240 equivalent representations, All have the same multiplication table, shown in Table \ref{tab:product} and the 16 representations of any primary act the same leave this primary invariant. So we can introduce $\GT_{i,O}, i\in\N_1^{30}$ that keeps $\Phi_{i,O}$ invariant along with the generated $\OB$. It also keeps all $\Phi_{i,j}$ within this primary resulting in $\Phi_{i,j'}$, $j,j'\in\N_1^{128}$. 

The remainder of this section derives all 21 terms of $\GT$, of which 14 are independent, using the notation of Bryant \cite{Bryant}. This will be compared to the Cartan formulation to uncover that the $\GT$ terms derived from the same remainder in (\ref{eqn:rho}), are automorphisms of the corresponding pseudo-octonion algebra, so may be considered as a sub-symmetry of $\GT$. 

Bryant defines $\GT$ as the Lie algebra with conjugations leaving the calibration 3-form $\Phi_{1,64}$ invariant, \cite{Bryant}, whereas \cite{HarveyLawson} used $\Phi_{1,23}$ and both represent octonions. Bryant then provides the same $\GT$ generators but states that he could not find an elegant proof and resorted to an explicit case. The generators in \cite{Bryant} can be converted to orthogonal matrices which can be converted to Bryant's form and represented in $\GA{7}$ as
\eqnb \begin{tabular}{lll} 
$A=\frac12(\e{23}-\e{45})$, &$H=\frac12(\e{45}-\e{67})$, &$A+H=\frac12(\e{23}-\e{67})$,\\[3pt]
$B=\frac12(-\e{13}-\e{46})$, &$I=\frac12(\e{46}+\e{57})$, &$B+I=\frac12(-\e{13}+\e{57})$,\\[3pt]
$C=\frac12(\e{12}+\e{47})$, &$J=\frac12(-\e{47}+\e{56})$, &$C+J=\frac12(\e{12}+\e{56})$,\\[3pt]
$D=\frac12(\e{15}+\e{26})$, &$K=\frac12(-\e{26}-\e{37})$, &$D+K=\frac12(-\e{15}-\e{37})$,\\[3pt]
$E=\frac12(\e{14}-\e{27})$, &$L=\frac12(-\e{27}+\e{36})$, &$E-L=\frac12(\e{14}-\e{36})$,\\[3pt]
$F=\frac12(\e{17}+\e{24})$, &$M=\frac12(-\e{17}+\e{35})$, &$F+M=\frac12(\e{24}+\e{35})$,\\[3pt]
$G=\frac12(-\e{16}-\e{25})$, &$N=\frac12(\e{25}-\e{34})$, &$G+N=\frac12(-\e{16}-\e{34})$.\\
\end{tabular} \label{eqn:G2}\eqne

\tabh\caption{Anti-symmetric Product Table for $\GT$}
\label{tab:product}    \centering \setlength{\tabcolsep}{0pt}
\begin{tabular}{|c|c|c|c|c|c|c|c|c|c|c|c|c|c|} \hline   
{\bf [,]} &{\bf  B} &{\bf C} &{\bf D} &{\bf E} &{\bf\;F\;} &{\bf G} &{\bf\;H\;} &{\bf\;I\;} &{\bf\;J\;} &{\bf K} &{\bf L} &{\bf M} &{\bf N}\\ \hline
{\bf A} &\sm\CJB &\sn\BI &\sm\LE &\sm\DKB &\sn N &\sn\FM &\sn0 &\sn J &\sm I &\sn L &\sm K &\sn N &\sm 2\FMB\\\hline
{\bf B} &\sn0 &\sm\AHB &\sm\FMB &\sm\GN &\sm K &\sm\LEB &\sm J &\sn0 &\sn H &\sn F &\sn\GN &\sn\DK &\sm\LEB\\ \hline
{\bf C} &. &\sn0 & \sm G &\sm 2F &\sn2E &\sn D &\sn I &\sm H &\sn0 &\sn\GN &\sm F &\sm E &\sm\DKB\\ \hline
{\bf D} &. &. &\sn0 &\sm H &\sn I &\sm 2\CJB &\sn E &\sm F &\sn G &\sn0 &\sm\AHB &\sm\BIB &\sn\CJ\\ \hline
{\bf E} &. &. &. &\sn0 &\sm 2C\; &\sn I &\sm D &\sm G &\sm F &\sm\AHB &\sn0 &\sn C &\sm\BIB\\ \hline
{\bf F} &. &. &. &. &\sn0 &\sn H &\sm G &\sn D \sn &\sn E &\sm B &\sn C &\sn0 &\sn A\\ \hline
{\bf G} &. &. &. &. &. &\sn0 &\sn F &\sn E &\sm D &\sm\CJB &\sm\BIB &\sn\AH &\sn0\\ \hline
{\bf H} &. &. &. &. &. &. &\sn0 &\sm 2J &\sn2I &\sn L &\sm K &\sm\GNB &\sn\FM\\ \hline
{\bf I} &. &. &. &. &. &. &. &\sn0 &\sm 2H\; &\sm\FMB &\sm N &\sn\DK &\sn L\\ \hline
{\bf J} &. &. &. &. &. &. &. &. &\sn0 &\sm N &\sn\FM &\sm\LEB &\sn K\\ \hline
{\bf K} &. &. &. &. &. &. &. &. &. &\sn0 &\sn2\AHB &\sm\BI &\sm J\\ \hline
{\bf L} &. &. &. &. &. &. &. &. &. &. &\sn0 &\sn\CJ &\sm I\\ \hline
{\bf M} &. &. &. &. &. &. &. &. &. &. &. &\sn0 &\sn A\\ \hline
{\bf N} &. &. &. &. &. &. &. &. &. &. &. &. &\sn0\\ \hline
\end{tabular}\tabe

The first two columns provide the basis for $\GT$ and last column contains dependent terms which provides the additional seven elements to cover all 21 rotations in $\Spin(7)$. Applying the pairs of terms in this representation as $90\degs$ rotations to $\Phi_{1,64}$ leaves it invariant. This is not the case for other $\Phi_O$ representations and $\GT$ must be transformed using the same transforms for $\Phi_{1,64}$. The Lie Product Table for $\GT$ is shown in Table \ref{tab:product} for the first 14 elements with only the upper half shown since the table is anti-symmetric.

These $\GT$ elements can be identified, alongside the representation from $(\ref{eqn:G2})$, from the enabling algebra given by the terms of $\Phi_{1,64}\hod$, as shown in Table \ref{tab:relation}. The three $90\degs$ rotation pairs, $R_{jklm}$, derived from the primary 4-forms, $\Phi_i\hod$, provide all 21 candidates for $\GT_{1,64}$. The relative signs for the $\GT$ elements, whether the two terms have the same sign or opposite signs, are provided by the parity in (\ref{eqn:auto}).

\tabb\caption{$\GT$ Relationship to $\Phi_1\hod$} 
\label{tab:relation} \centering
{\renewcommand{\arraystretch}{1.1} \small
\begin{tabular}{|@{}c@{}|c|c|c|}  \hline  
{$\bf\:\Phi_{1,64}\hod$\ term\:} &{\bf Normal Rotations} &{\bf Mixed Rotations} &{\bf Outer Rotations}\\ \hline
$-\e{1247}$ &$C=\frac12(\e{12}+\e{47})$ &$E=\frac12(\e{14}-\e{27})$ &$F=\frac12(\e{17}+\e{24})$\\[1pt] \hline
$-\e{1256}$ &$C+J=\frac12(\e{12}+\e{56})$ &$-D=\frac12(\e{15}-\e{26})$ &$-G=\frac12(\e{16}+\e{25})$\\[1pt] \hline
$-\e{1346}$ &$-B=\frac12(\e{13}+\e{46})$ &$E-L=\frac12(\e{14}-\e{36})$ &$-G-N=\frac12(\e{16}+\e{34})$\\[1pt] \hline
$\e{1357}$ &$-B-I=\frac12(\e{13}-\e{57})$ &$-D-K=\frac12(\e{15}+\e{37})$ &$-M=\frac12(\e{17}-\e{35})$\\[1pt] \hline
$\e{2345}$ &$A=\frac12(\e{23}-\e{45})$ &$F+M=\frac12(\e{24}+\e{35})$ &$N=\frac12(\e{25}-\e{34})$\\[1pt] \hline
$\e{2367}$ &$A+H=\frac12(\e{23}-\e{67})$ &$-K=\frac12(\e{26}+\e{37})$ &$-L=\frac12(\e{27}-\e{36})$\\[1pt] \hline
$\e{4567}$ &$H=\frac12(\e{45}-\e{67})$ &$I=\frac12(\e{46}+\e{57})$ &$-J=\frac12(\e{47}-\e{56})$\\[1pt] \hline
\end{tabular}}\tabe

\begin{figure}[ht] \begin{center}
  \begin{tabular}{c}
  \includegraphics{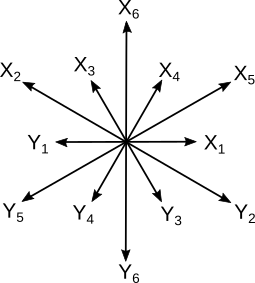}
  \end{tabular}
  \caption{\bf Cartan Root Diagram for $\GT$} \label{fig:cartan}
\end{center} \end{figure}

Table \ref{tab:relation} for any $\Phi_{i,j}$ has at least one invariant row. The remainder from the Classification Theorem leaves one row invariant because $(\ref{eqn:rho})$ is invariant for the three rotations involving the remainder's dual. It is only when the remainder in (\ref{eqn:rho}) is zero that all rows are invariant. So the six algebras, $\PB4, \PB8, \PB{10}, \PB{12}, \PB{14}, \PB{16}$, provide three automorphisms each that are subsets of $\GT$ and leave $\Phi$ invariant. This can be partially related to the 6 fold symmetry of the Cartan root system diagram for $\GT$, \cite{CanHowe}, shown in Figure \ref{fig:cartan} using construction and destruction matrices, $X_i$ and $Y_i, i\in\N_1^6$, respectively, derived from the $g$ matrices in Humphreys \cite{Humphreys}. The mapping is

\eqb \arraycolsep=1.4pt\begin{array}{lllllll} X_1 = g_2^T, &X_2 = g_{1,-2}, &X_3 = -g_1^T, &X_4 = g_3, &X_5 = -g_{2,-3}, &X_6 = g_{1,-3} \\
Y_1 = g_2, &Y_2 = g_{1,-2}^T, &Y_3 = -g_1, &Y_4 = g_3^T, &Y_5 = -g_{2,-3}^T, &Y_6 = g_{1,-3} \\
H_1 = [X_1 , Y_1], &&H_2 = [X_2 , Y_2].&&&&\end{array} \eqe
These generate the $\GT$ multiplication table in \cite{CanHowe} replicated in Table \ref{tab:cartan}.
Matching these $90\degs$ rotations to $\Phi_{14}$ leads to the following symmetry identification
\eqb E-L = X_1-Y_1, \quid D+K = X_2-Y_2, \quid H = Y_3-X_3. \eqe
Hence the corresponding row in Table \ref{tab:relation} that is left invariant by the pseudo-octonion algebra that is given by the remainder in the first column matches one of the 6-fold axes in Figure \ref{fig:cartan}. Thus the symmetry of octonions is broken and it is only when all rows in the table are invariant that the full invariance is achieved. This is a partial mapping because other $\GT$ terms can not be represented and three terms only match half of the symmetry provided by the $\GT$ Cartan root diagram. This could relate to the groups of six primaries representing different pseudo-octonion algebras.

The identification of $\GT$ in Table \ref{tab:relation} can be replaced with a construction of $\GT$ due to an inherent asymmetry in $\GT$. Any $\Phi_O$ representation generates a table similar to Table \ref{tab:relation} except the order of the rows and terms may differ. The one invariant is that the row that contains M always occurs in the row identified from the Classification Theorem as the term which changes sign when starting from the primary ($\Phi_{11}$ and $\Phi_{20}$ are also easily recognised). This observation allows for the construction of the representation of the $\GT$ algebra, up to the overall sign. For $\Phi_1$, the negated term is $\e{246}$ and $-\e{1234567}\e{246} = \e{1357}$ and this is consistent across all octonion representations in $\Phi_1$. 

\tabb\caption{$\GT$ Multiplication Structure in Cartan representation}
\label{tab:cartan}\footnotesize   \centering \setlength{\tabcolsep}{2pt}
\begin{tabular}{|c|c|c|c|c|c|c|c|c|c|c|c|c|c|} \hline
{\bf [,]} &\;{\bf \HH2}\; &{\bf \XX1} &{\bf \YY1} &{\bf \XX2} &{\bf \YY2} &{\bf \XX3} &{\bf \YY3} &{\bf \XX4} &{\bf \YY4} &{\bf \XX5} &{\bf \YY5} &{\bf \XX6} &{\bf \YY6}\\ \hline
{\bf \HH1}  &\sn0 &\sn2\XX1 &\sm2\YY1 &\sm3\XX2 &\sn3\YY2 &\sm\XX3 &\sn\YY3 &\sn\XX4 &\sm\YY4 &\sn3\XX5 &\sm3\YY5 &\sn0 &\sn0\\ \hline
{\bf \HH2}  &\sn0 &\sm\XX1 &\sn\YY1 &\sn2\XX2 &\sm2\YY2 &\sn\XX3 &\sm\YY3 &\sn0 &\sn0 &\sm\XX5 &\sn\YY5 &\sn\XX6 &\sm \YY6\\ \hline
{\bf \XX1}  &. &\sn0 &\sn\HH1 &\sn\XX3 &\sn0 &\sn2\XX4 &\sm3\YY2 &\sm3\XX5 &\sm2\YY3 &\sn0 &\sn\YY4 &\sn0 &\sn0\\ \hline
{\bf \YY1}  &. &. &\sn0 &\sn0 &\sm\YY3 &\sn3\XX2 &\sm2\YY4 &\sn2\XX3 &\sn3\YY5 &\sm\XX4 &\sn0 &\sn0 &\sn0\\ \hline
{\bf \XX2}  &. &. &. &\sn0 &\sn\HH2 &\sn0 &\sn\YY1 &\sn0 &\sn0 &\sm\XX6 &\sn0 &\sn0 &\sn\YY5\\ \hline
{\bf \YY2}  &. &. &. &. &\sn0 &\sm\XX1 &\sn0 &\sn0 &\sn0 &\sn0 &\sn\YY6 &\sm \XX5 &\sn0\\ \hline
{\bf \XX3}  &. &. &. &. &. &\sn0 &\sn\HH1{+}3\HH2 &\sm3\XX6 &\sn2\YY1 &\sn0 &\sn0 &\sn0 &\sn\YY4\\ \hline
{\bf \YY3}  &. &. &. &. &. &. &\sn0 &\sm2\XX1 &\sn3\YY6 &\sn0 &\sn0 &\sm\XX4 &\sn0\\ \hline
{\bf \XX4}  &. &. &. &. &. &. &. &\sn0 &\sn\HH1{+}2\HH2 &\sn0 &\sm\YY1 &\sn0 &\sm\YY3\\ \hline
{\bf \YY4}  &. &. &. &. &. &. &. &. &\sn0 &\sn\XX1 &\sn0 &\sn\XX3 &\sn0\\ \hline
{\bf \XX5}  &. &. &. &. &. &. &. &. &. &\sn0 &\sn\HH1{+}\HH2 &\sn0 &\sm\YY2\\ \hline
{\bf \YY5}  &. &. &. &. &. &. &. &. &. &. &\sn0 &\sn\XX2 &\sn0\\ \hline
{\bf \XX6}  &. &. &. &. &. &. &. &. &. &. &. &\sn0 &\sn\HH1{+}2\HH2\\ \hline
{\bf \YY6}  &. &. &. &. &. &. &. &. &. &. &. &. &\sn0\\ \hline
\end{tabular}\tabe

The $M$ row is identified as above but it is uncertain which of the three rotations provides $M$ so the three candidates are trialed. First take the Lie product with all other terms then take cross products of these results. Using $Z_i, i\in\N_1^9$ to designate $\GT$ terms in the columns of Table \ref{tab:G2}, then for rotation elements $Z_1, Z_2$ and $Z_3$, the Lie $M$ products are $[M, Z_i] = Z_{i+3}, i\in\N_1^3$, and cross terms are $[Z_4, Z_5] = Z_7, [Z_4, Z_6] = Z_8, [Z_5, Z_6] = Z_9$. These are shown for $\Phi_{1,64}$ in Table \ref{tab:G2} and a pattern is obvious. Three rows have two zeros in the cross products, one which is the $M$ row, and the other four rows have repeated $M$ Product pairs, $B+I$ and $D+K$. These are the other elements in the $M$ row which indicates a symmetry in $\GT$ that could be hard to separate. Fortunately, of the two pairs, the $C+J$ and $E-L$ terms always occur first in the term ordering shown in the table.

\tabb\caption{Construction of $\GT$}  
\label{tab:G2}\small    \centering
\begin{tabular}{|c|c|c|c|} \hline    
{\bf$\Phi_{1,64}\hod$\ term} &{\bf Enabling Rotations} &{\bf M Products} &{\bf Cross Products}\\ \hline
$\e{1247}$ & $C, E, F$ & $ -E, C, 0$ &$ -2F, 0, 0$\\ \hline
$\e{1256}$ & $C+J, -D, -G$ & $-L, B+I, -A-H$ &$ G, 2K, -C-J$\\ \hline
$\e{1346}$ &$-B, E-L, -G-N$ & $-D-K, -J, +H$ &$ G+N, E-L, -2I$\\ \hline
$-\e{1357}$ &  $-B-I, -D-K, -M$ & $-2(D+K), 2(B+I), 0$ &$ -8M, 0, 0$\\ \hline
$-\e{2345}$ &  $A, F+M, N$ & $N, 0, A$ &$ 0, -2(F+M), 0$\\ \hline
$-\e{2367}$ & $A+H, -K, -L$ & $-G, B+I, -C-J$ &$ -L, -2D, A+H$\\ \hline
$-\e{4567}$ & $H, I , -J$ & $-G-N, D+K, L-E$ &$ -J, 2B, H$\\ \hline
\end{tabular}\tabe

Given any $\Phi_i\hod$, the relationship table can be formed and the $M$ row for each trial group identified. Of the three trials trying to find which column represents $M$, the three rows with two zeros in the cross terms can be distinguished from the other four containing the $M$ Product pairs. Excluding the $8M$ row from the two zeros rows leaves the $C/E$ and $A/N$ rows. The $M$ Products here are asymmetric with $[M, C] = -E$ while all other products are positive. Hence $C$ and $E$ can be identified as well as $F$ from the Cross Products column. Because the $A$ and $N$ products are symmetric, these are left ambiguous for now but $F+M$ can be identified from the cross products. 

The 2-form terms within $C$ and $E$ can now isolate the first rows with pairs $B+I$ and $D+K$ within the $M$ Products, respectively. The other row of the pair has these same 2-form terms in $J$ and $L$ but as previously mentioned the $C$ and $E$ terms always appear first in the relationship table. With the location of the paired elements and the $C$ and $E$ terms identified then $D, G, L$ and $A+H$ can be isolated from the first row of the $B+I$ pair and $B, J, H$ from the first $D+K$ pair. The paired rows provide $K$ and $I$ using these paired locations. Finally, $A+H$ and $H$ allows the ambiguity of $A$ and $N$ to be settled. This process generates one of three solutions with $M$ sometimes replaced by $D+K$ or $B+I$. But the $M$ solution always occurs so automorphisms are needed to distinguish the outer rotation solution as being the correct one. This construction is messy but is included to highlight that $\Cl(7)$ provides a complete derivation of the Lie algebra $\GT$.

\section{Summary}
The calibrations that define unique cross products in seven dimensions are usually expressed using differential geometry. Mapping the exterior algebra to the basis of $\GA{7}$ allows the cross product and Hodge-star operations to become products in $\GA{}$. It also allows a form to be defined that is invertible but not related to the orthogonal group. This form provides an easy method to uncover all representations of the octonions since the form $\rho_O^2 = -1$ in (\ref{eqn:rho}) is invariant under all $\Pin(7)$ transformations. Furthermore, this form can be used to analyse all sign combinations $\rho_{i,j}$ not in $\rho_O$ which have squares with a remainder from the enabling algebra. This distinguishes a number of algebras in groups of six primaries and over all primaries six unique algebras can be identified, categorised by the number of non-associative triple products, called $\PB4, \PB8, \PB{10}, \PB{12}, \PB{14}, \PB{16}$. These were called the pseudo-octonion algebras and were shown to break the symmetry of octonions. The construction also shows that there is no possibility of an associative algebra being derived from a 7-D cross product 3-form. The pseudo-octonion algebras each have 12 zero divisors which makes them akin to sedenions and appropriate for a naming scheme, $\PB{k}$, where $k$ is the number of unique, non-associative ordered triples. The sedenions and higher level Cayley-Dickson algebras up to level 8 can be denoted as $\SB=\PB{252}$, $\PB{3612}$, $\PB{36204}$, $\PB{319788}$, and $\PB{2678508}$, which are all power-associative algebras. The split octonions, ${\tilde\OB}$, which are also power-associative, can be interpreted as the seventh algebra in the series and designated as $\PB{28}$. 

Geometric algebra has derived the octonion and pseudo-octonion algebras by identifying the terms of $\Phi_i\hod$ as a subalgebra of $\Spin(7)$. The $\Phi_i\hod$ automorphisms keep $\rho_{i,j}^2$ from $(\ref{eqn:rho})$ invariant for the remainder term's dual, or for all terms if the remainder is zero. These terms identify elements of $\GT$ and the algebra can be constructed in a systematic way using Lie bracket operations. In this respect the connection between octonions, pseudo-octonions and $\GT$ is intrenched in $\Spin(7)$ since 4-form rotations separate into $\GT$ rotation terms and other terms that cancel. The coassociative 4-form, $\Phi_{i,O}\hod$, generates 16 representations of $\GT$ for any primary which only differ by sign changes for the $\GT$ terms. The connection between $\Spin(7)$, octonions and $\GT$ is usually summarised as the dimension count relationship $21 = 7 +14$. The definition of $\GT$ can be expressed more precisely as $\GT_{i,j}$ for any $\Phi_{i,j}\in\Phi_O$, and these consist of the automorphisms of $\Spin(7)$ that keep $\Phi_{i,j}$ invariant but other representations within the $i$ primary change between themselves. This can be extended to $\GT_{i,j}$ as a subalgebra of $\GT_i$ provided by the remainder row in Table \ref{tab:relation} which are the automorphisms that preserve $\Phi_{i,j}$. But the dimensional relationship breaks down in this case.

The connection between $\Spin(7)$ and octonions can be extended to $\Spin(15)$ and sedenions and $100$ pseudo-sedenion algebras. This scheme should extend to higher dimensions, so that $\Spin(31)$ would generate the second level of sedenions, $\PB{3612}$, and further subalgebras. As expected, there is no $\GT$ equivalent for $\SB$ and presumably the series and it is hard to imagine a generalisation since the dual of the 3-form is a 12-form in $\GA{15}$ of which the terms are not closed under multiplication.

This work was verified with the use of geometric algebra and quaternion, octonion, sedenion, et cetera calculators written in Python. Quaternions were used to verify the Clifford calculator using Euler angles and these generalised to arbitrary dimensions and compared to matrix expansions of rotations for up to seven dimensions. The octonion calculator also has 15 dimensions but can practically only support up to 8, giving sedenions up to eight levels, called $\PB{2678508}$ above. It also allows for unitary basis elements, $\u{i}^2 = +1$, which generates split octonions. The github URL for the calculators is \hfil\break \centerline{\url{https://github.com/GPWilmot/geoalg}.}

An advantage of geometric algebra is the conciseness of the notation which has allowed all 480 representations of the octonions and $\GT$ to be enumerated. Matrix representations for $R_{ij}$ were not needed and many of the calculations can be done on the fly using the simple product and swapping rules. Another advantage is the geometric interpretation of each operation and the geometric understanding provided. For example, the subalgebra that $\Phi\hod$ spans does not involve 6-forms because a set of three rotations is inconsistent with octonion multiplication or for the six related pseudo-octonion algebras. I would like to acknowledge the invaluable assistance and encouragement of Professor Derek Abbott with the publishing of this article. I also gratefully acknowledge the support of an Australian Government Research Training Program Scholarship.

\end{document}